\documentclass[11pt,oneside,reqno]{article}
\usepackage[utf8]{inputenc}
\usepackage{braket}
\usepackage{slashed}
\usepackage{graphicx}
\usepackage{amsmath}
\usepackage{amsfonts} 
\usepackage{marginnote}
\usepackage{xcolor}
\usepackage{hyperref}
\usepackage{multicol}
\usepackage[margin=1.25in]{geometry}
\hypersetup{
    colorlinks=true,
    linkcolor=blue,
    filecolor=magenta,      
    urlcolor=cyan,
    citecolor=red,
    pdftitle={symbolicmethod},
    pdfpagemode=FullScreen,
    }
\graphicspath{ {./images/} }

\numberwithin{equation}{section}
\title{An Operational Calculus Generalization of Ramanujan's Master Theorem}
\author{Zachary P. Bradshaw$^{1}$ and Christophe Vignat$^{1,2}$}
\date{
    $^1$Tulane University, Department of Mathematics, New Orleans, LA, USA\\
    $^2$LSS, CentraleSup\'elec, Universit\'e Paris Saclay, France\\[2ex]
}

\begin{document}

\maketitle
\begin{abstract}
    We give a formal extension of Ramanujan's master theorem using operational methods. The resulting identity transforms the computation of a product of integrals on the half-line to the computation of a Laplace transform. Since the identity is purely formal, we show consistency of this operational approach with various standard calculus results, followed by several examples to illustrate the power of the extension. We then briefly discuss the connection between Ramanujan's master theorem and identities of Hardy and Carr before extending the latter identities in the same way we extended Ramanujan's. Finally, we generalize our results, producing additional interesting identities as a corollary.
\end{abstract}
\section{Introduction}
A surprising formal justification of Ramanujan's master theorem \cite{amdeberhan2012,berndt2012,hardy1959}, a remarkably simple tool for evaluating the loop integrals arising from Feynman diagrams \cite{MoBFeynman}, can be given using operational calculus methods \cite{roman2019}; we aim to extend this theorem symbolically. Indeed, such a justification was given in \cite{gorska2012} and we include it here for completeness. Let $\phi_n=(-1)^n/n!$ and suppose
\begin{align}
    f(x)=\sum_{n=0}^\infty\phi_ng(n)x^n,
\end{align}
where $g(n)$ is a sequence that we assume has a natural extension as a function $z\to g(z)$ defined over the complex plane. Symbolically, we may write
\begin{align}
    f(x)&=\sum_{n=0}^\infty\phi_nx^ne^{n\partial}g(0)\\
    &=e^{-xe^\partial}g(0),
\end{align}
where the symbol $\partial$ denotes a derivative in the variable associated to the extension of $g$, so that $e^{n\partial}$ becomes the translation operator and we have the operational rule $e^{n\partial} g(0)=g(n)$. Let us now substitute this result into the Mellin transform of $f$. We have
\begin{align}
    \mathcal{M}(f)(s)&:=\int_0^\infty x^{s-1}f(x)\ dx\\
    &=\int_0^\infty x^{s-1}e^{-xe^\partial}g(0)\ dx\\
    &=\int_0^\infty x^{s-1}e^{-xe^\partial}\ dx\ g(0).\label{eq:changeofvariable}
\end{align}
Now performing the substitution $x\to xe^{-\partial}$, we have
\begin{align}
    \mathcal{M}(f)(s)&=\int_0^\infty x^{s-1} e^{-x} e^{-s\partial}\ dx\ g(0)\\
    &=\int_0^\infty x^{s-1} e^{-x}\ dx\ g(-s)\\
    &=\Gamma(s)g(-s),
\end{align}
and this is precisely Ramanujan's master theorem, which was proven rigorously under natural conditions by Hardy in \cite{hardy1959}. 

We emphasize here that there are many steps in this calculation which are not rigorous, so that this operational approach is nothing more than a mathematical curiosity. However, we claim that this operational method can be modified to include functions other than $x^{s-1}$, and the resulting formal identity is a helpful tool in the study of definite integrals, including the evaluation of Feynman diagrams. Indeed, consider the integral
\begin{align}
    I=\int_0^\infty h(x)f(x)\ dx,
\end{align}
for which the Laplace transform $H$ of $h$ exists and $f$ is as before. Then using the operational rule $e^{n\partial}g(0)=g(n)$,
\begin{align}
    I=\int_0^\infty h(x) e^{-xe^\partial}\ dx\ g(0),
\end{align}
and now we recognize the integral as the Laplace transform $H(s)$ of $h$ evaluated at $e^\partial$ and applied to $g(0)$. That is,
\begin{align}\label{eq:Lrmt}
    \int_0^\infty h(x)f(x)\ dx=H(e^\partial)g(0).
\end{align}

We recover Ramanujan's master theorem by setting $h(x)=x^{s-1}$. To see this, observe that the Laplace transform of $h$ becomes $H(p)=\Gamma(s)/p^s$, so that $H(e^\partial)g(0)=\Gamma(s)e^{-s\partial}g(0)=\Gamma(s)g(-s)$. An advantage to this approach in deriving Ramanujan's master theorem is that the operator valued change of variables needed in \eqref{eq:changeofvariable} is no longer needed.

It is straightforward to extend \eqref{eq:Lrmt} to multivariate integrals, as the construction is similar to the univariate case. Consider the function
\begin{align}
    f(x_1,\ldots,x_k)=\sum_{n_1,\ldots,n_k=0}^\infty \phi_{n_1\cdots n_k}g(n_1,\ldots,n_k)x_1^{n_1}\cdots x_k^{n_k},
\end{align}
where we have made the identification $\phi_{n_1\cdots n_k}=\phi_{n_1}\cdots\phi_{n_k}$. Let $\partial_i$ denote the derivative in the $i$-th component of $g$. Then we may write
\begin{align}
    g(n_1,\ldots,n_k)=e^{n_1\partial_1+\cdots+n_k\partial_k}g(0,\ldots,0),
\end{align}
so that the function $f$ becomes
\begin{align}
    f(x_1,\ldots,x_k)=e^{-x_1e^{\partial_1}}\cdots e^{-x_ke^{\partial_k}}g(0,\ldots,0).
\end{align}

Now plugging this result into the multivariate version of the integral in \eqref{eq:Lrmt}, we obtain
\begin{align}
    \int_{\mathbb{R}_+^k}h(\Vec{x})f(\Vec{x})\ dx^k&=\int_{\mathbb{R}_+^k}h(\Vec{x})e^{-x_1e^{\partial_1}}\cdots e^{-x_ke^{\partial_k}}\ dx^k\ g(0,\ldots,0)\\
    &=H(e^{\partial_1},\ldots,e^{\partial_k})g(0,\ldots,0)\label{eq:multiLrmt},
\end{align}
where $H(s_1,\ldots,s_k)$ denotes the $k$-dimensional Laplace transform of $h$. As a special case, we recover a multivariate version of Ramanujan's master theorem \cite{amdeberhan2012} when we set $h(x_1,\ldots,x_k)=x_1^{s_1-1}\cdots x_k^{s_k-1}$. Indeed, in this case we have
\begin{align}
    \int_{\mathbb{R}_+^k}x_1^{s_1-1}\cdots x_k^{s_k-1}f(\Vec{x})\ dx^k&=\Gamma(s_1)\cdots\Gamma(s_k)e^{-s_1\partial_1}\cdots e^{-s_k\partial_k}g(0,\ldots,0)\\
    &=\Gamma(s_1)\cdots\Gamma(s_k)g(-s_1,\ldots,-s_k),
\end{align}
for which setting $k=1$ reproduces Ramanujan's master theorem.

Although the above argument is purely formal, we claim that the resulting identity is a useful tool in the study of definite integrals. This is because it easily produces formal identities for integrands in which $h$ has a sufficiently nice Laplace transform. In such cases, $H(s)$ can be expanded as a power series, and we recover a possible series representation for the definite integral in question. However, the validity of the result has to be checked by other means.

The rest of this paper is dedicated to the study of the formal identity \eqref{eq:Lrmt}. In Section~\ref{sec:consistency}, we show that this identity is consistent with standard calculus results such as integration by parts, the inverse Laplace transform, the fundamental theorem of calculus, and a change of variables. In Section~\ref{sec:examples} we give several examples that illustrate the power of this method, including the evaluation of an integral arising from the bubble diagram. In Section~\ref{sec:HRC} we connect Ramanujan's master theorem with identities due to Hardy and Carr and extend them in the same way we have extended the master theorem. We discuss analogs of \eqref{eq:Lrmt} with transforms other than the Laplace transform in Section~\ref{sec:generalization}. Finally, in Section~\ref{sec:conclusion} we give concluding remarks.

\section{Consistency with Standard Calculus Results}\label{sec:consistency}
Without a rigorous justification of the formal identity \eqref{eq:Lrmt}, it is necessary to verify its consistency with standard results from calculus, and this is the purpose of this section. We will show that this operational method is consistent with a change of variables, integration by parts, the fundamental theorem of calculus, and the inverse Laplace transform.
\subsection{Change of Variables} Let us consider the integral
\begin{align}
    I=\int_0^\infty h(x)f(x)\ dx,
\end{align}
where $f$ again has the form $f(x)=\sum_n\phi_ng(n)x^n$. Assume there is a parametrization $x(t)$ such that $x(0)=0$ and $x(t)\to\infty$ as $t\to\infty$. By performing a change of variables $x\to x(t)$, this integral becomes
\begin{align}
    \int_0^\infty h(x)f(x)\ dx=\int_0^\infty h(x(t))f(x(t))x'(t)\ dt.
\end{align}
We will verify that the formal identity \eqref{eq:Lrmt} produces consistent results on the left and right hand sides. On the left hand side, \eqref{eq:Lrmt} produces
\begin{align}
    \int_0^\infty h(x)f(x)\ dx=H(e^\partial)g(0).
\end{align}
Meanwhile, on the right hand side, we have
\begin{align}
    \int_0^\infty h(x(t))f(x(t))x'(t)\ dt=\int_0^\infty h(x(t)) e^{-x(t)e^\partial}x'(t)\ dt\ g(0),
\end{align}
and reverting to the original variable by the substitution $x(t)\to x$, we have
\begin{align}
    \int_0^\infty h(x(t))f(x(t))x'(t)\ dt&=\int_0^\infty h(x)e^{-xe^\partial}\ dx\ g(0)\\
    &=H(e^\partial)g(0).
\end{align}
Thus, the operational method is consistent with a change of variables.

\subsection{Integration by Parts}
Suppose that $f$ and $h$ are such that $h(0)=0$ and the boundary term vanishes when we perform an integration by parts. That is,
\begin{align}
    \int_0^\infty h'(x)f(x)\ dx=-\int_0^\infty h(x)f'(x)\ dx.
\end{align}
We will check for the consistency of the operational method with this identity. On the left hand side, we see that the Laplace transform of $h'(x)$ must be computed. Recall that we have $\mathcal{L}\{h'(x)\}(s)=sH(s)-h(0)$ where $H(s)$ is the Laplace transform of $h(x)$. Then by our formal identity \eqref{eq:Lrmt}, the left hand side evaluates to
\begin{align}
    (e^\partial H(e^\partial)-h(0))g(0)=e^\partial H(e^\partial)g(0)-h(0)g(0),
\end{align}
but we have made the assumption that $h(0)=0$ so that we are left with $e^\partial H(e^\partial)g(0)$. Meanwhile, in the right hand integral, we observe that \begin{align}f'(x)&=\sum_{n=1}^\infty\phi_ng(n)nx^{n-1}\\
&=\sum_{n=0}^\infty\phi_{n+1}(n+1)g(n+1)x^n\\
&=-\sum_{n=0}^\infty\phi_ng(n+1)x^n\\
&=\sum_{n=0}^\infty\phi_n\tilde g(n) x^n,
\end{align}
where we have made the identification $\tilde g(n):=-g(n+1)$. Then the right hand integral is given by \eqref{eq:Lrmt} and evaluates to
\begin{align}
    -H(e^\partial)\tilde g(0)=H(e^\partial)g(1)=e^\partial H(e^\partial) g(0).
\end{align}
Thus, under the conditions we have outlined, the operational identity \eqref{eq:Lrmt} is consistent with integration by parts.

\subsection{Fundamental Theorem of Calculus} In this section we show consistency of the fundamental theorem of calculus with a generalized version of \eqref{eq:Lrmt} to any interval of integration $(a,b)$. Indeed, let us show that
\begin{align}
    \int_a^b \frac{d}{dx}(h(x)f(x))\ dx=h(b)f(b)-h(a)f(a),
\end{align}
where $f$ again has the form $f(x)=\sum_n\phi_ng(n)x^n$. Starting with the integral, we apply the product rule
\begin{align}
    \int_a^b \frac{d}{dx}(h(x)f(x))\ dx=\int_a^b h'(x)f(x)\ dx+\int_a^b h(x)f'(x)\ dx.
\end{align}
Now let us restrict the integration interval in the Laplace transform to the interval $(a,b)$, so that
\begin{align}
    H(p)=\int_a^b e^{-px}h(x)\ dx.
\end{align}
Then
\begin{align}
    \int_a^b h'(x)f(x)\ dx&=\int_a^b h'(x)e^{-xe^\partial}\ dx\ g(0)\\
    &=\bigg([e^{-xe^\partial}h(x)]_a^b+e^\partial\int_a^be^{-xe^\partial}h(x)\bigg)g(0)\\
    &=\bigg(e^{-be^\partial}h(b)-e^{-ae^\partial}h(a)+e^\partial H(e^\partial)\bigg)g(0)\\
    &=f(b)h(b)-f(a)h(a)+H(e^\partial)g(1)
\end{align}
Meanwhile, $f'(x)=-\sum_n\phi_ng(n+1)x^n=-e^\partial e^{-xe^\partial}g(0)$, so that
\begin{align}
    \int_a^b h(x)f'(x)\ dx=-H(e^\partial)g(1).
\end{align}
Putting this together, we have
\begin{align}
    \int_a^b\frac{d}{dx}(h(x)f(x))\ dx&=h(b)f(b)-h(a)f(a)+H(e^\partial)g(1)-H(e^\partial)g(1)\\
    &=h(b)f(b)-h(a)f(a).
\end{align}
Thus, we have established consistency with the fundamental theorem of calculus. Note that setting $a=0$ and $b=\infty$ recovers the case of consistency with \eqref{eq:Lrmt}.

\subsection{Inverse Laplace Transform}
Suppose $F(s)$ is the Laplace transform of some function $f(x)$, and we wish to take the inverse Laplace transform to recover $f(x)$. There is a well-known formula for this purpose due to Emil L. Post \cite{post1930}:
\begin{align}\label{eq:inverseL}
    f(x)=\lim_{m\to\infty}\phi_m\bigg(\frac{m}{x}\bigg)^{m+1}F^{(m)}\bigg(\frac{m}{x}\bigg).
\end{align}
Here we will check that this formula is consistent with our formal identity \eqref{eq:Lrmt}. Let $h(x)=e^{-sx}$. Then we have
\begin{align}
    F(s)=\int_0^\infty h(x)f(x)\ dx=H(e^\partial)g(0).
\end{align}
Now, $H(p)=\frac{1}{p+s}$, so that
\begin{align}
    F(s)=\frac{1}{s+e^\partial}g(0).
\end{align}
Taking the $m$-th derivative with respect to $s$ yields
\begin{align}
    F^{(m)}(s)=\frac{(-1)^mm!}{(s+e^\partial)^{m+1}}g(0).
\end{align}
Then
\begin{align}
    \lim_{m\to\infty}\phi_m\bigg(\frac{m}{x}\bigg)^{m+1}F^{(m)}\bigg(\frac{m}{x}\bigg)&=\lim_{m\to\infty}\frac{(-1)^m}{m!}\bigg(\frac{m}{x}\bigg)^{m+1}\frac{(-1)^mm!}{(\frac{m}{x}+e^\partial)^{m+1}}g(0)\\
    &=\lim_{m\to\infty}\frac{1}{(1+\frac{xe^\partial}{m})^{m+1}}g(0)\\
    &=e^{-xe^\partial}g(0)\\
    &=f(x).
\end{align}
Thus, \eqref{eq:inverseL} is consistent with our formal identity \eqref{eq:Lrmt}.

\section{Examples}\label{sec:examples}

It is the purpose of this section to show that there are many interesting examples which illustrate the power of the operational method we have outlined.

\subsection{Sum of Hermite Polynomials}
Let $H_n(x)$ denote the physicist's Hermite polynomials defined by the exponential generating function
\begin{align}
    e^{2xt-t^2}=\sum_{n=0}^\infty H_n(x)\frac{t^n}{n!}.
\end{align}
Then our operational method can be used to show that
\begin{align} \label{eq:hermitesum}
    \sum_{n=0}^\infty H_n(x)=\frac{\sqrt{\pi}}{2}e^{x^2-x+1/4}(1-\text{erf}(1/2-x)),
\end{align}
where erf$(x)$ denotes the error function. Surprisingly, the symbolic programming language, \texttt{Mathematica}, is unable to produce \eqref{eq:hermitesum}. Now, using (\cite{prudnikov1992}, Volume 4, 2.2.1.5), we have
\begin{align}
    \int_0^\infty e^{2xt-t^2}e^{-t} dt=\frac{\sqrt{\pi}}{2}e^{x^2-x+1/4}(1-\text{erf}(1/2-x)).
\end{align}
Note that there is a typo in the entry 2.2.1.5 of \cite{prudnikov1992}. The $\sqrt{\frac{\pi}{p}}$ should read $\sqrt{\frac{\pi}{a}}$. Meanwhile, letting $f(t)=e^{-t}$ and $h(t)=e^{2xt-t^2}=\sum_{n=0}^\infty H_n(x)\frac{t^n}{n!}$, and noting that $g(n)=1$ and the Laplace transform of $h$ is $\sum_{n=0}^\infty H_n(x)s^{-n-1}$, we apply \eqref{eq:Lrmt} to deduce that
\begin{align}
    \sum_{n=0}^\infty H_n(x)=\frac{\sqrt{\pi}}{2}e^{x^2-x+1/4}(1-\text{erf}(1/2-x)).
\end{align}

\subsection{Three Double Integral Examples}

\begin{enumerate}
\item In (\cite{prudnikov1986}, Volume 1, 3.1.3.5), the following two-dimensional Laplace transform is computed:
\begin{align}
    \int_0^\infty\int_0^\infty e^{-px-qy}\frac{1}{\sqrt{xy}}e^{-\frac{1}{x^2y}}\ dxdy=\pi\sqrt{\frac{2}{pq}}e^{-2p^{1/2}q^{1/4}}.
\end{align}
Let $f(x,y)=\sum_{m,n=0}^\infty\phi_{n,m}g(n,m)x^ny^m$. By \eqref{eq:multiLrmt}, we have
\begin{align}
    \int_0^\infty\int_0^\infty f(x,y)\frac{1}{\sqrt{xy}}e^{-\frac{1}{x^2y}}\ dxdy=\pi\sqrt{2}\sum_{n=0}^\infty\phi_n2^ng\bigg(\frac{n-1}{2},\frac{n-2}{4}\bigg).
\end{align}

\item Similarly, it is shown in (\cite{prudnikov1986}, Volume 1, 3.1.3.61) that
\begin{align}
    \int_0^\infty\int_0^\infty e^{-px-qy}\frac{\cos(2\sqrt{axy})}{\sqrt{xy}}\ dxdy=\frac{\pi}{\sqrt{a+pq}},
\end{align}
which we recognize as the two-dimensional Laplace transform of the function $\frac{\cos(2\sqrt{axy})}{\sqrt{xy}}$. Then by \eqref{eq:multiLrmt}, we have
\begin{align}
    \int_0^\infty\int_0^\infty f(x,y)\frac{\cos(2\sqrt{axy})}{\sqrt{xy}}\ dxdy&=\frac{\pi}{\sqrt{a}}\sum_{n=0}^\infty\frac{\Gamma(n+1/2)}{\Gamma(n+1)\Gamma(1/2)}\bigg(-\frac{1}{a}\bigg)^ng(n,n)\\
    &=\sqrt{\frac{\pi}{a}}\sum_{n=0}^\infty\frac{\Gamma(n+1/2)}{\Gamma(n+1)}\bigg(-\frac{1}{a}\bigg)^ng(n,n),
\end{align}
which depends only on the diagonal values of $g$. Moreover, after a change of variables, we have
\begin{align}
    4\int_0^\infty\int_0^\infty f(x^2,y^2)\cos(2\sqrt{a}xy)\ dxdy=\sqrt{\frac{\pi}{a}}\sum_{n=0}^\infty\frac{\Gamma(n+1/2)}{\Gamma(n+1)}\bigg(-\frac{1}{a}\bigg)^ng(n,n).
\end{align}
Letting $a=1$ and $f(x,y)=e^{-x-y}$, we have $g(n,n)=1$ and
\begin{align}
    4\int_0^\infty\int_0^\infty e^{-x^2-y^2}\cos(2xy)\ dxdy&=\sqrt{\pi}\sum_{n=0}^\infty\frac{\Gamma(n+1/2)}{\Gamma(n+1)}(-1)^n\\
    &=\sqrt{\pi}\sqrt{\frac{\pi}{2}}\\
    &=\frac{\pi}{\sqrt{2}},
\end{align}
which is consistent with the corresponding \texttt{Mathematica} calculation.

\item Another example we can consider comes from (\cite{prudnikov1986}, Volume 1, 3.1.6.1), which says
\begin{align}
    \int_0^\infty\int_0^\infty (2\gamma+\log(xy))e^{-px-qy}\ dxdy=-\frac{\log(pq)}{pq},
\end{align}
where $\gamma$ is Euler's constant. We recognize the integral as a double Laplace transform so that with $h(x)=(2\gamma+\log(xy))$, we have $H(e^\partial_1,e^\partial_2)=-(\partial_1+\partial_2)e^{-\partial_1-\partial_2}$. From this, we deduce that
\begin{align}
    \int_0^\infty\int_0^\infty (2\gamma+\log(xy))f(x,y)\ dxdy=-g^{(1,0)}(-1,-1)-g^{(0,1)}(-1,-1),
\end{align}
where $g^{(k,l)}=\partial_1^k\partial_2^lg$. Now let us consider the function \begin{align}
    f(x,y)&=\frac{1}{(1+x+y)^3}\\
    &=\frac12\sum_{m,n}\phi_{m,n}x^my^n(m+n+2)!\ ,
\end{align}
so that $g(m,n)=\frac12\Gamma(m+n+3)$, and we have
\begin{align}
    \int_0^\infty\int_0^\infty \frac{(2\gamma+\log(xy))}{(1+x+y)^3}\ dxdy&=-g^{(1,0)}(-1,-1)-g^{(0,1)}(-1,-1)\\
    &=-\frac12\Gamma'(1)-\frac12\Gamma'(1)\\
    &=\gamma.
\end{align}
More generally, with $p\ge3$, define
\begin{align}
    f(x,y)&=\frac{1}{(1+x+y)^p}\\
    &=\frac{1}{(p-1)!}\sum_{m,n}\phi_{m,n}x^my^n(m+n+p-1)!.
\end{align}
Then
\begin{align}
    g(m,n)=\frac{1}{(p-1)!}\Gamma(m+n+p)=\frac{\Gamma(m+n+p)}{\Gamma(p)},
\end{align}
from which we deduce
\begin{align}
    \int_0^\infty\int_0^\infty \frac{(2\gamma+\log(xy))}{(1+x+y)^p}\ dxdy&=-g^{(1,0)}(-1,-1)-g^{(0,1)}(-1,-1)\\
    &=-\frac{2\Gamma'(p-2)}{\Gamma(p)}\\
    &=-\frac{2\psi(p-2)}{(p-1)(p-2)},
\end{align}
where $\psi$ denotes the digamma function.

\end{enumerate}

\subsection{A Triple Integral Example} From (\cite{prudnikov1986}, Volume 1, 3.2.3.2), we have
\begin{align}
    \int_0^\infty\int_0^\infty\int_0^\infty\frac{(xyz)^{\frac{\nu-2}{2}}}{(xy+xz+yz)^{\frac{\nu+1}{2}}}e^{-px-qy-rz}\ dxdydz=\frac{8\pi\Gamma(\nu/2+1)}{\nu(\nu-1)(\nu-2)}(\sqrt{p}+\sqrt{q}+\sqrt{r})^{2-\nu},
\end{align}
and we recognize the integral as a triple Laplace transform. Now, using the Schwinger parametrization, we have
\begin{align}
    (\sqrt{p}+\sqrt{q}+\sqrt{r})^{2-\nu}=\frac{1}{\Gamma(\nu-2)}\int_0^\infty x^{\nu-3}e^{-x(\sqrt{p}+\sqrt{q}+\sqrt{r})}\ dx.
\end{align}
Let $f(x,y,z)=\sum_{m,n,k}\phi_{m,n,k}g(m,n,k)x^{m}y^{n}z^{k}=e^{-xe^{\partial_1}-ye^{\partial_2}-ze^{\partial_3}}g(0,0,0)$. Then we have
\begin{align}
    \int_0^\infty\int_0^\infty&\int_0^\infty\frac{(xyz)^{\frac{\nu-2}{2}}}{(xy+xz+yz)^{\frac{\nu+1}{2}}}f(x,y,z)\ dxdydz\\
    &=\frac{8\pi\Gamma(\nu/2+1)}{\nu(\nu-1)(\nu-2)}\frac{1}{\Gamma(\nu-2)}\int_0^\infty x^{\nu-3}e^{-x(e^{\partial_1/2}+e^{\partial_2/2}+e^{\partial_3/2})}\ dx\ g(0,0,0). \label{eq:parametrizedform}
\end{align}
Define a new function by $\Bar{f}(x,y,z)=\sum_{m,n,k=0}^\infty\phi_{m,n,k}g(m/2,n/2,k/2)x^my^nz^k$. Then \eqref{eq:parametrizedform} becomes
\begin{align}
    \frac{8\pi\Gamma(\nu/2+1)}{\Gamma(\nu+1)}\int_0^\infty x^{\nu-3}\bar{f}(x,x,x)\ dx.
\end{align}
To illustrate the use of this formula, consider the normalized confluent hypergeometric $U$-function (\cite{NIST:DLMF}, entry 13.4.4)
\begin{align}
    &\frac{1}{2^p(x+y+z)^{p/2}}U\bigg(\frac{p}{2},\frac12;\frac{1}{4(x+y+z)}\bigg)\\
    &=\frac{1}{\Gamma(p)}\int_0^\infty\lambda^{p-1}e^{-\lambda}e^{-\lambda^2(x+y+z)}\ d\lambda\\
    &=\frac{1}{\Gamma(p)}\int_0^\infty e^{-\lambda}\sum_{l,m,n}\phi_{l,m,n}x^ly^mz^n\lambda^{2l+2m+2n+p-1}\ d\lambda\\
    &=\sum_{l,m,n}\phi_{l,m,n}x^ly^mz^n\frac{\Gamma(2l+2m+2n+p)}{\Gamma(p)},
\end{align}
so that $g(l,m,n)=\frac{\Gamma(2l+2m+2n+p)}{\Gamma(p)}$ and
\begin{align}
    \bar{f}(x,y,z)&=\sum_{m,n,l}\phi_{m,n,l}x^my^nz^l\frac{\Gamma(l+m+n+p/2)}{\Gamma(p)}\\
    &=\frac{1}{(1+x+y+z)^{p/2}}.
\end{align}
Then from
\begin{align}
    \int_0^\infty x^{\nu-3}\bar{f}(x,x,x)\ dx&=\int_0^\infty x^{\nu-3}\frac{1}{(1+3x)^{p/2}}\ dx\\
    &=3^{2-\nu}\frac{\Gamma(\nu-2)\Gamma(2+\frac{p}{2}-\nu)}{\Gamma(\frac{p}{2})},
\end{align}
which holds for $2<\nu<\frac{p+4}{2}$, $p>0$, we deduce that
\begin{align}
    \int_0^\infty\int_0^\infty\int_0^\infty& \frac{(xyz)^{\frac{\nu-2}{2}}}{(xy+xz+yz)^{\frac{\nu+1}{2}}}\frac{1}{2^p(x+y+z)^{p/2}}U\bigg(\frac{p}{2},\frac12;\frac{1}{4(x+y+z)}\bigg)\ dxdydz\\
    &=\frac{8\pi\Gamma(\nu/2+1)}{\Gamma(\nu+1)}3^{2-\nu}\frac{\Gamma(\nu-2)\Gamma(2+\frac{p}{2}-\nu)}{\Gamma(\frac{p}{2})}\\
    &=3^{2-\nu}\frac{8\pi\Gamma(\nu/2+1)}{\Gamma(\nu+1)}B(\nu-2,2+p/2-\nu),
\end{align}
where $B(x,y)$ denotes Euler's beta function.

\subsection{Modified Bessel Functions}
Let $I_\nu$ denote the modified Bessel functions of the first kind \cite{abramowitz1972}. It can be shown \cite{prudnikov1992} that they satisfy the interesting Laplace transform identity
\begin{align}
    \mathcal{L}\bigg\{\frac{\pi}{2}I_0\bigg(\frac{ax}{2}\bigg)^2\bigg\}(p)=\frac{1}{p}K\bigg(\frac{a}{p}\bigg),
\end{align}
where $K(x)$ is the complete elliptic function
\begin{align}
    K(x)=\frac{\pi}{2}\sum_{n=0}^\infty \binom{2n}{n}^2\bigg(\frac{x}{4}\bigg)^{2n}.
\end{align}
Let $f(x)=\sum_{n=0}^\infty\phi_ng(n)x^n$ as before and apply \eqref{eq:Lrmt} to deduce
\begin{align}
    \int_0^\infty f(x)\bigg(I_0\bigg(\frac{ax}{2}\bigg)\bigg)^2\ dx=\sum_{n=0}^\infty\binom{2n}{n}^2\bigg(\frac{a}{4}\bigg)^{2n}g(-2n-1).
\end{align}
As another example, we note that \cite{prudnikov1992}
\begin{align}
    \mathcal{L}\bigg\{\frac{\pi a}{2}I_0\bigg(\frac{ax}{2}\bigg)I_1\bigg(\frac{ax}{2}\bigg)\bigg\}(p)=K\bigg(\frac{a}{p}\bigg)-\frac{\pi}{2},
\end{align}
so that applying \eqref{eq:Lrmt} produces the identity
\begin{align}
    \int_0^\infty f(x)I_0\bigg(\frac{ax}{2}\bigg)I_1\bigg(\frac{ax}{2}\bigg)\ dx=\frac{1}{a}\sum_{n=1}^\infty\binom{2n}{n}^2\bigg(\frac{a}{4}\bigg)^{2n}g(-2n).
\end{align}

\subsection{Dedekind Eta Function}
Let $\eta(x)$ denote the Dedekind $\eta$-function. It is shown in \cite{glasser2009} that the Laplace transform of this function with an imaginary argument is
\begin{align}
    \mathcal{L}\{\eta(ix)\}(t)=\sqrt{\frac{\pi}{t}}\frac{\sinh(2\sqrt{\frac{\pi t}{3}})}{\cosh(\sqrt{3\pi t})},
\end{align}
and we have the series expansion
\begin{align}
    \frac{1}{x}\frac{\sinh(\alpha x)}{\cosh(\beta x)}=\frac{1}{2}\sum_{n=0}^\infty\frac{E_n\bigg(\frac{\beta-\alpha}{2\beta}\bigg)-E_n\bigg(\frac{\beta+\alpha}{2\beta}\bigg)}{n!}(-2\beta x)^n,
\end{align}
where $E_n$ denotes the $n$-th Euler polynomials defined by the exponential generating function
\begin{align}
    \frac{2e^{xz}}{e^x+1}=\sum_{n=0}^\infty\frac{E_n(z)}{n!}x^n.
\end{align}
Since it can be easily checked that $E_{2n}(1/6)-E_{2n}(5/6)=0$, we have
\begin{align}
    \mathcal{L}\{\eta(ix)\}(t)=\sqrt{\pi}\sum_{n=0}^\infty\frac{E_{2n+1}(1/6)-E_{2n+1}(5/6)}{(2n+1)!}(-2\sqrt{3\pi})^{2n+1}t^n.
\end{align}
Thus, by applying \eqref{eq:Lrmt}, we deduce
\begin{align}
    \int_0^\infty\eta(ix)f(x)\ dx=\sqrt{\pi}\sum_{n=0}^\infty\frac{E_{2n+1}(5/6)-E_{2n+1}(1/6)}{(2n+1)!}(2\sqrt{3\pi})^{2n+1}g(n),
\end{align}
where we have again assumed that $f(x)=\sum_{n=0}^\infty\phi_ng(n)x^n$.

\subsection{Massless Bubble Diagram Integral} 
\begin{figure}[ptb]
\begin{center}
\includegraphics[
width=4.0in
]{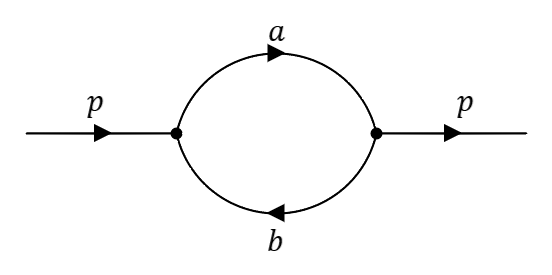}
\end{center}
\caption{The Bubble Diagram}
\label{fig:bubble}%
\end{figure}

Consider the massless bubble diagram in Figure~\ref{fig:bubble}. The Schwinger parametrization of the associated integral is \begin{align}
    I&=\frac{(-1)^{-D/2}}{\Gamma(a)\Gamma(b)}\int_0^\infty\int_0^\infty\frac{x^{a-1}y^{b-1}}{(x+y)^{D/2}}e^{-p^2\frac{xy}{x+y}}\ dxdy\\
    &=\frac{(-1)^{-D/2}}{\Gamma(a)\Gamma(b)}\int_0^\infty\int_0^\infty\frac{x^{a-D/2-1}y^{b-1}}{(1-(-y/x))^{D/2}}e^{-\frac{(-p^2x)(-y/x)}{1-(-y/x)}}\ dxdy\\
    &=\frac{(-1)^{-D/2}}{\Gamma(a)\Gamma(b)}\int_0^\infty\int_0^\infty x^{a-D/2-1}y^{b-1}\sum_{n=0}^\infty L_n^{(D/2-1)}(-p^2 x)\bigg(-\frac{y}{x}\bigg)^n\ dxdy \label{eq:pluginto},
\end{align}
where $L_n^{(\alpha)}$ denotes the generalized Laguerre polynomial defined by the generating function
\begin{align}
    \frac{e^{-\frac{tx}{1-t}}}{(1-t)^{\alpha+1}}=\sum_{n=0}^\infty L_n^{(\alpha)}(x) t^n.
\end{align}
It can be shown that the closed form of the generalized Laguerre polynomial is 
\begin{align}
    L_n^{(\alpha)}(x)&=\sum_{k=0}^n\phi_k\binom{n+\alpha}{n-k}x^k\\
    &=\sum_{k=0}^\infty\phi_k\frac{\Gamma(n+\alpha+1)}{\Gamma(n-k+1)\Gamma(n+\alpha+k+1)}x^n.
\end{align}
Plugging this in to \eqref{eq:pluginto}, we have
\begin{align}
    I=\frac{(-1)^{-D/2}}{\Gamma(a)\Gamma(b)}\int_0^\infty\int_0^\infty x^{a-D/2-1}y^{b-1}\sum_{n,k=0}^\infty\phi_{n,k}\frac{(-p^2)^k\Gamma(n+D/2)\Gamma(n+1)}{\Gamma(n-k+1)\Gamma(k+D/2)}x^{k-n}y^ndxdy.
\end{align}
Now performing the substitution $(u,v)=(y/x,x)$, we have $dudv=\frac{1}{x}dxdy$, $x=v$, and $y=uv$, so that
\begin{align}
    I&=\frac{(-1)^{-D/2}}{\Gamma(a)\Gamma(b)}\int_0^\infty\int_0^\infty u^{b-1}v^{a+b-D/2-1}\sum_{n,k=0}^\infty\phi_{n,k}\frac{(-p^2)^k\Gamma(n+D/2)\Gamma(n+1)}{\Gamma(n-k+1)\Gamma(k+D/2)}u^nv^kdudv.
\end{align}
In order to apply \eqref{eq:multiLrmt}, we let $h(u,v)=u^{b-1}v^{a+b-D/2-1}$, so that its Laplace transform is $H(s,t)=\Gamma(b)\Gamma(a+b-c)s^{-b}t^{c-a-b}$. Then with 
\begin{align}
    g(n,k)=(-1)^k(p^2)^k\frac{\Gamma(n+D/2)\Gamma(n+1)}{\Gamma(n-k+1)\Gamma(k+D/2)},
\end{align}
we have by \eqref{eq:multiLrmt} that
\begin{align}
    I&=\frac{(-1)^{-D/2}}{\Gamma(a)\Gamma(b)}\Gamma(b)\Gamma(a+b-c)g(-b,c-a-b)\\
    &=(-1)^{-a-b}(p^2)^{D/2-a-b}\frac{\Gamma(a+b-D/2)\Gamma(D/2-b)\Gamma(1-b)}{\Gamma(a)\Gamma(a-D/2+1)\Gamma(D-a-b)},
\end{align}
which, using the gamma function identity $\Gamma(z-n)=(-1)^{n-1}\frac{\Gamma(-z)\Gamma(1+z)}{\Gamma(n+1-z)}$, $n\in\mathbb{Z}$, transforms to
\begin{align}
    I=(-1)^{-D/2}(p^2)^{D/2-a-b}\frac{\Gamma(a+b-D/2)\Gamma(D/2-a)\Gamma(D/2-b)}{\Gamma(a)\Gamma(b)\Gamma(D-a-b)},
\end{align}
and this is the correct evaluation of the integral. Thus, the operational approach can be used to evaluate Feynman diagrams.

\section{On Identities of Hardy, Ramanujan, and Carr}\label{sec:HRC} We have already generalized Ramanujan's master theorem,
\begin{align}
    \int_0^\infty x^{\nu-1}\sum_{n=0}^\infty\phi_ng(n)x^n\ dx=\Gamma(\nu)g(-\nu)
\end{align}
to the formal identity \eqref{eq:Lrmt}. Hardy proved a different version of this theorem, which can be obtained from Ramanujan's version by making the identification $g(n)=n!\tilde g(n)$. Indeed, this produces
\begin{align}
    \int_0^\infty x^{\nu-1}\sum_{n=0}^\infty(-1)^n\tilde g(n)x^n\ dx&=\Gamma(\nu)\Gamma(1-\nu)\tilde g(-\nu)\label{eq:hardyformula}\\
    &=\frac{\pi}{\sin(\pi\nu)}\tilde g(-\nu),\label{eq:hardy}
\end{align}
where in the last equality we have applied Euler's reflection formula.

\subsection{Hardy's Formula}
What happens if we extend Hardy's formula as we did Ramanujan's? Replacing $x^{\nu-1}$ by an arbitrary function $h$, we have
\begin{align}
    \int_0^\infty h(x)\sum_{n=0}^\infty(-1)^n\tilde g(n)x^n\ dx&=\int_0^\infty h(x)\sum_{n=0}^\infty(-1)^nx^ne^{n\partial}\ dx\ \tilde g(0)\\
    &=\int_0^\infty h(x)\frac{1}{1+xe^\partial}\ dx\ \tilde g(0)\\
    &=\pi e^\partial Hi(-e^{-\partial})\tilde g(0),\label{eq:hilbert}
\end{align}
where 
\begin{align}
    Hi(p)=\frac{1}{\pi}\int_0^\infty \frac{h(x)}{p-x}\ dx
\end{align} 
is the one-sided Hilbert transform of $h$. Alternatively, we can plug $g(n)=n!\tilde g(n)$ directly into \eqref{eq:Lrmt} to obtain
\begin{align}
    \int_0^\infty h(x)\sum_{n=0}^\infty(-1)^n\tilde g(n)x^n\ dx=H(e^\partial)[\Gamma(1)\tilde g(0)],\label{eq:hilbertlaplace}
\end{align}
where the square bracket is to indicate that the operator $\partial$ applies to the entirety of $\Gamma(1)\tilde g(0)$. For example, $e^{k\partial}[\Gamma(1)\tilde g(0)]=\Gamma(k+1)\tilde g(k)$. Moreover, from \eqref{eq:hilbert} and \eqref{eq:hilbertlaplace}, we deduce that
\begin{align}
    Hi(-e^{-\partial})\tilde g(0)=\frac{1}{\pi}e^{-\partial}H(e^\partial)[\Gamma(1)\tilde g(0)].
\end{align}

\subsection{Carr's Identity} Suppose that we further alter Ramanujan's master theorem by ridding the power series in Hardy's formula \eqref{eq:hardyformula} of the factor of $(-1)^n$. We can accomplish this by making the identification $\tilde g(n)=\sin(\frac{(2n+1)\pi}{2})\hat{g}(n)$. Then \eqref{eq:hardyformula} becomes
\begin{align}
    \int_0^\infty x^{\nu-1}\sum_{n=0}^\infty\hat{g}(n)x^n\ dx&=\frac{\pi}{\sin(\pi\nu)}\sin\bigg(\frac{(1-2\nu)\pi}{2}\bigg)\hat{g}(-\nu)\\
    &=\pi\cot(\pi\nu)\hat{g}(-\nu),
\end{align}
which is a variation of Carr's identity (\cite{carr1886}, entry 2709),
\begin{align}
    \int_0^\infty x^{\nu-1}e^{-x}\sum_{k=0}^\infty g(k)x^k\ dx=\Gamma(\nu)\sum_{k=0}^\infty (\nu)_kg(k),
\end{align}
where $(x)_n$ denotes the Pochhammer symbol \cite{roman2019}. Again we replace $x^{\nu-1}$ with an arbitrary function $h(x)$ to find
\begin{align}
    \int_0^\infty h(x)\sum_{n=0}^\infty \hat{g}(n)x^n\ dx&=\int_0^\infty h(x)\sum_{n=0}^\infty x^ne^{n\partial}\ dx\ \hat{g}(0)\\
    &=\int_0^\infty h(x)\frac{1}{1-xe^\partial}\ dx\ \hat{g}(0)\\
    &=\pi e^\partial Hi(e^{-\partial})\hat{g}(0).
\end{align}

\section{Generalization to Other Integral Transforms}\label{sec:generalization}
The examples in Section~\ref{sec:examples} made use of the formal identities \eqref{eq:Lrmt} and \eqref{eq:multiLrmt}. However, as we saw in Section~\ref{sec:HRC}, the operational method outlined here is capable of handling many different choices for the form of the function $f(x)$. In many cases, a similar result is produced, but the Laplace transform is replaced by another integral transform. We will first give two more examples that illustrate the power of the operational method, and then move to a more general setting.

\subsection{Cosine Transform}
Let us consider a function $f$ with the expansion
\begin{align}
    f(x)=\sum_{n=0}^\infty\phi_{2n}(ix)^{2n}g(2n).
\end{align}
It is shown by O. Atale in \cite{atale2022} that the Mellin transform of $f$ is given by
\begin{align} \label{eq:atale}
    \int_0^\infty x^{s-1}f(x)\ dx=\Gamma(s)g(-s)\cos\bigg(\frac{\pi s}{2}\bigg),
\end{align}
but we will re-derive this result by constructing a cosine transform analog of \eqref{eq:Lrmt}. Observe that
\begin{align}
    f(x)=\sum_{n=0}^\infty\phi_{2n}(ix)^{2n}e^{2n\partial}g(0)=\cos(xe^\partial)g(0).
\end{align}
Then
\begin{align}
    \int_0^\infty h(x)f(x)\ dx&=\int_0^\infty h(x)\cos(xe^\partial)\ dx\ g(0)\\
    &=H_c(e^\partial)g(0),
\end{align}
where $H_c$ denotes the (half) cosine transform of $h$. To see that this reproduces \eqref{eq:atale}, let $h(x)=x^{s-1}$. Then $H_c(p)=\frac{\Gamma(s)}{p^s}\cos(\frac{\pi s}{2})$, from which we deduce
\begin{align}
    \int_0^\infty x^{s-1}f(x)\ dx&=\Gamma(s)\cos\bigg(\frac{\pi s}{2}\bigg)e^{-s\partial}g(0)\\
    &=\Gamma(s)g(-s)\cos\bigg(\frac{\pi s}{2}\bigg).
\end{align}
\subsection{Zeta Function Identity}
There is another interesting result due to Atale in \cite{atale2021}. Let $f(x)=\sum_{n=0}^\infty\phi_ng(n)x^n$ and define $F(x)=\sum_{m=1}^\infty f(mx)$. Then Atale shows that the Mellin transform of $F(x)$ is given by
\begin{align}\label{eq:atale2}
    \int_0^\infty x^{s-1}F(x)\ dx=\zeta(s)\Gamma(s)g(-s).
\end{align}
Here we will use the operational method we have outlined to generalize this identity. Writing $f(x)=e^{-xe^\partial}g(0)$, we observe that $F(x)=\sum_{m=1}^\infty e^{-mxe^\partial}g(0)$. Then
\begin{align}
    \int_0^\infty h(x)F(x)\ dx&=\sum_{m=1}^\infty\int_0^\infty e^{-mxe^\partial}h(x)\ dx\ g(0)\\
    &=\sum_{m=1}^\infty H(me^\partial)g(0).
\end{align}
Notice that the choice $h(x)=x^{s-1}$ recovers \eqref{eq:atale2}.

\subsection{General Setting}
An interesting question is whether the previous example can be generalized to other transforms. Indeed, let $f(x)=\sum_{n=0}^\infty\phi_ng(n)x^n$ as before and assume that $\psi(s)$ satisfies the hypothesis of Ramanujan's master theorem. Define a new function by $f_\psi(x)=\sum_n\phi_ng(n)\psi(n)x^n$ and define the kernel of a transform by
\begin{align}
    K_\psi(x,y)=\sum_n\phi_n(xy)^n\psi(n).
\end{align}
Then we have
\begin{align}
    f_\psi(x)=K_\psi(x,e^\partial)g(0).
\end{align}
\sloppy Denote the transform of $h(x)$ with respect to the kernel $K_\psi(x,y)$ by $H_\psi(y):=\int_0^\infty h(x)K_\psi(x,y)\ dx$. Define $F_\psi(x)=\sum_m f_\psi(mx)$. Then
\begin{align}
    \int_0^\infty h(x)F_\psi(x)\ dx&=\sum_m\int_0^\infty h(x)K_\psi(mx,e^\partial)\ dx\ g(0)\\
    &=\sum_m\int_0^\infty h(x)K_\psi(x,me^\partial)\ dx\ g(0)\\
    &=\sum_m H_\psi(me^\partial)g(0).
\end{align}
Consider the particular case $h(x)=x^{s-1}$. We have
\begin{align}
    H(y)&=\int_0^\infty x^{s-1}\sum_n\phi_n\psi(n)y^nx^n\ dx\\
    &=\psi(-s)y^{-s}\Gamma(s),
\end{align}
where in the last equality we have applied Ramanujan's master theorem.
Then
\begin{align}
    \int_0^\infty x^{s-1}F_\psi(x)\ dx&=\sum_m\psi(-s)(me^\partial)^{-s}\Gamma(s)g(0)\\
    &=\psi(-s)\Gamma(s)g(-s)\zeta(s).
\end{align}

Let us consider a few examples.

\begin{enumerate}
\item Let $\psi(s)=1$. Then we recover the result of Atale \cite{atale2021}:
\begin{align}
    \int_0^\infty x^{s-1}F_\psi(x)\ dx=\Gamma(s)g(-s)\zeta(s).
\end{align}

\item Let $\psi(s)=\cos(\frac{s\pi}{2})$ and notice that $\psi(2n)=(-1)^n$ while $\psi(2n+1)=0$. We should therefore recover the case involving the cosine transform. Indeed, we have
\begin{align}
    \int_0^\infty x^{s-1}F_\psi(x)\ dx=\Gamma(s)g(-s)\cos\bigg(\frac{s\pi}{2}\bigg)\zeta(s).
\end{align}
Now applying Riemann's functional equation, we have 
\begin{align}
    \int_0^\infty x^{s-1}F_\psi(x)\ dx=\frac{(2\pi)^s\zeta(1-s)}{2}g(-s).
\end{align}

\item Similarly, let $\psi(s)=\sin(\frac{s\pi}{2})$ and notice that $\psi(2n)=0$ while $\psi(2n+1)=(-1)^n$. We should therefore recover the case involving the sine transform. Indeed, we have
\begin{align}
    \int_0^\infty x^{s-1}F_\psi(x)\ dx=-\Gamma(s)g(-s)\sin\bigg(\frac{s\pi}{2}\bigg)\zeta(s).
\end{align}

\item Let \begin{align}\psi(s)=-\frac{2}{\sqrt{\pi}}\sin(\pi s/2)\frac{\Gamma(1+s/2)}{\Gamma(1/2+s/2)}=-2\frac{\sin(\pi s/2)}{B(1/2,s/2)},
\end{align} 
where $B(x,y)$ denotes the Euler beta function, and notice that $\psi(2n)=0$ while $\psi(2n+1)=-\frac{(-1)^n\Gamma(2n+2)}{(\Gamma(n+1))^24^n}$. Then
\begin{align}
    K_\psi(x,y)=\sum_{n=0}^\infty\phi_n\frac{(yx)^{2n+1}}{\Gamma(n+1)4^n}=J_0(yx)yx,
\end{align}
which is the kernel of the product of the Hankel transform with $\nu=0$ and a factor of $y$. In this case, we find
\begin{align}
    \int_0^\infty x^{s-1}F_\psi(x)\ dx&=\bigg(\frac{\pi}{2}\bigg)^s\frac{\zeta(1-s)}{B(1/2,s/2)}g(-s).
\end{align}
\end{enumerate}

\section{Conclusion}\label{sec:conclusion}
It is clear that the operational approach outlined in this work is capable of producing interesting identities. The caveat is that the calculations we have done are purely formal, so that any identity derived here would need to be checked for validity by other means. Therefore, a natural question is whether any of the formal results discussed here can be made rigorous. Still, this approach can be used as a powerful tool for computing potential representations of a given definite integral.

There is a method related to Ramanujan's master theorem that was originally constructed for use in the computation of definite integrals arising from the evaluation of Feynman diagrams \cite{Klauber,Sakurai} called the method of brackets \cite{part1,part2,MoBFeynman}. An interesting question for future work is whether \eqref{eq:Lrmt} and \eqref{eq:multiLrmt} can be used to extend this method. However, the current formulation is already simple, and may prove itself to be a valuable tool in the study of high energy particle physics without the added simplicity of the framework of the method of brackets.

\newpage
\bibliographystyle{unsrt}
\bibliography{Ref}

\end{document}